\documentclass{amsart}
\usepackage[T1]{fontenc}
\usepackage[latin1]{inputenc}
\usepackage{color}
\usepackage[enableskew,vcentermath]{youngtab}
\usepackage{amsfonts}
\usepackage{hhline}
\usepackage{amsthm}
\usepackage{amsmath}

\begin{document}

\title[Multiplicity-Free Skew Characters]{On Multiplicity-Free Skew Characters And The Schubert Calculus}
\author[C. Gutschwager]{Christian Gutschwager}
\address{Institut für Algebra, Zahlentheorie und Diskrete Mathematik, Leibniz Universität Hannover,  Welfengarten 1, D-30167 Hannover}
\email{gutschwager (at) math (dot) uni-hannover (dot) de}

\newtheorem{Le}{Lemma}[section]
\newtheorem{Ko}[Le]{Lemma}
\newtheorem{Sa}[Le]{Theorem}
\newtheorem{pro}[Le]{Proposition}
\newtheorem{Bem}[Le]{Remark}
\newtheorem{Def}[Le]{Definition}
\newtheorem{Bsp}[Le]{Example}


\begin{abstract} In this paper we introduce a partial order on the set of skew characters of the symmetric group which we use to classify the multiplicity-free skew characters. Furthermore we give a short and easy proof that the Schubert calculus is equivalent to that of skew characters in the following sense: If we decompose the product of two Schubert classes we get the same as if we decompose a skew character and replace the irreducible characters by Schubert classes of the `inverse' partitions (Theorem~\ref{Sa:Mainskewschub}).
\end{abstract}
\maketitle

\section{Introduction}
In this paper we introduce a partial order on the set of skew characters by proving an inequality of the Littlewood-Richardson coefficients (Theorem~\ref{Sa:schieb}). This we use to classify the multiplicity-free skew characters of the symmetric group $S_n$ (Theorem~\ref{Sa:Mainskew}), i.e. skew characters $[\lambda/\mu]$ for which in the decomposition $[\lambda/\mu]=\sum_\nu c(\lambda;\mu,\nu) [\nu]$ with irreducible characters $[\nu]$ all coefficients are $0$ or $1$.

By comparing the LR-coefficients of skew characters $[\lambda/\mu]=\sum_\nu c(\lambda;\mu,\nu) [\nu]$  and skew Schur functions, which are given by $s_{\lambda/\mu}=\sum_\nu c(\lambda;\mu,\nu) s_{\nu}$, we see that skew characters and skew Schur functions are equivalent.

Furthermore, as pointed out below, the classification of multiplicity-free skew characters is by Theorem~\ref{Sa:Mainskewschub} equivalent to the classification of multiplicity-free products of Schubert classes which was done by Thomas and Yong in \cite{Tho}. The difference between the proofs in \cite{Tho} and this paper is that we don't use ad hoc constructions of LR-fillings to prove multiplicity in the cases where multiplicity appears but rather reduce the problem by Theorem~\ref{Sa:schieb} to a few small cases. However by Remark~\ref{Rem}.(\ref{Rem3}) we still can produce easily the LR-fillings which give multiplicity in the cases where multiplicity appears.

Skew characters whose skew diagram $\mathcal{D}$ decompose into disconnected diagrams $\mathcal{A},\mathcal{B}$ are equivalent to the product of the  characters of the disconnected diagrams induced to a larger symmetric group and equivalent to the product of Schur functions. For the Schur functions this gives: $s_\mathcal{D}=s_\mathcal{A}\cdot s_\mathcal{B}$. The classification of  mul\-ti\-plicity-free products of Schur functions was done by Stembridge in \cite{Ste}. However we will not give a new proof of Stembridge's classification but rather use his classification in our proofs.

In Section~\ref{sec:Schubert} we will give a short and easy proof for the equivalence of skew characters and the Schubert calculus, which is a product of two Schubert classes indexed by partitions $\mu,\nu$ defined by $\sigma_\mu \cdot \sigma_\nu=\sum_{\lambda\subseteq (k^l)} c(\lambda;\mu,\nu) \sigma_\lambda$ with some positive integers $k,l$. Theorem~\ref{Sa:Mainskewschub} tells us that this sum is the skew character $[((k^l)/\nu)^\circ/\mu]$ if we replace in the sum the Schubert classes $\sigma_\lambda$ by characters $[((k^l)/\lambda)^\circ]$. The classification when this product is multiplicity-free is given by Thomas and Yong in \cite{Tho} and is equivalent to our classification of multiplicity-free skew characters (Theorem~\ref{Sa:Mainskew}).

Finally we will arrive at Theorem~\ref{Sa:Mainskewshort} which gives a classification of multiplicity-free skew characters  which looks like Stembridge's classification for multiplicity-free products of Schur functions.

\section{Notation and Littlewood-Richardson-Symmetries}
We mostly follow the standard notation in \cite{Sag}. A partition $\lambda=(\lambda_1,\lambda_2,\ldots,\lambda_l)$ is a weakly decreasing sequence of non-negative integers. For the length we write $l(\lambda)=l$ and the sum $\left|\lambda\right|=\sum_i \lambda_i$. With a partition $\lambda$ we associate a diagram, which we also denote by $\lambda$, containing $\lambda_i$ left-justified boxes in the $i$-th row and we use matrix-style coordinates to refer to the boxes.

The conjugate $\lambda'$ of $\lambda$ is the diagram which has $\lambda_i$ boxes in the $i$-th column. By $\lambda \cup (n)$ we refer to the partition $(\lambda_1,\ldots,\lambda_i,n,\lambda_{i+1},\ldots)$ when $\lambda_i\geq n \geq \lambda_{i+1}$, i.e. inserting a row with $n$ boxes into the diagram $\lambda$. $\lambda + (1^n)$ is the partition $(\lambda_1+1,\lambda_2+1,\ldots,\lambda_n+1,\lambda_{n+1},\ldots)$, i.e. inserting a column with $n$ boxes into the diagram $\lambda$. Both operations can be generalized to partitions $\nu$ instead of $(n)$ (resp. $(1^n)$) and one gets $\left(\lambda \cup \nu\right)'=\lambda'+\nu'$.

For $\mu \subseteq \lambda$ we define the skew diagram $\lambda/\mu$ as the difference of the diagrams $\lambda$ and $\mu$ defined as the difference of the sets of boxes. Rotation of $\lambda/\mu$ by $180^\circ$ yields a skew diagram $(\lambda/\mu)^\circ$ which is well defined up to translation. A skew tableau $T$ is a skew diagram in which the boxes are replaced by positive integers.  We refer to the entry in box $(i,j)$ as $T(i,j)$. A semistandard tableau of shape $\lambda/\mu$ is a filling of $\lambda/\mu$ with positive integers such that the following expressions hold for all $(i,j)$ for which they are defined: $T(i,j)<T(i+1,j)$ and $T(i,j)\leq T(i,j+1)$. The content of a semistandard tableau $T$ is $\nu=(\nu_1,\ldots)$ if the number of occurrences of the entry $i$ in $T$ is $\nu_i$. The reverse row word of a tableau $T$ is the sequence obtained by reading the entries of $T$ from right to left and top to bottom starting at the first row. Such a sequence is said to be a lattice word if for all $i,n \geq1$ the number of occurrences of $i$ among the first $n$ terms is at least the number of occurrences of $i+1$ among these terms. The Littlewood-Richardson (LR-) coefficient $c(\lambda;\mu,\nu)$ equals the number of semistandard tableaux of shape $\lambda/\mu$ with content $\nu$ such that the reverse row word is a lattice word. We will call those tableaux LR-tableaux. The LR-coefficients play an important role in different contexts (see \cite{Sag}).

The irreducible characters $[\lambda]$ of the symmetric group $S_n$ belong to Young's natural representation of the Specht module $S^\lambda$, with partitions $\lambda\vdash n$. The skew character $[\lambda/\mu]$ to a skew diagram $\lambda/\mu$ is defined by the LR-coefficients:
\[ [\lambda/\mu]=\sum_\nu c(\lambda;\mu,\nu) [\nu] \]

There are many known symmetries of the LR-coefficients (see \cite{Sag}).

We have that $c(\lambda;\mu,\nu)=c(\lambda;\nu,\mu)$. If the skew diagrams of $\lambda/\mu$ and $\alpha/\beta$ are the same up to translation then we get $c(\lambda;\mu,\nu)=c(\alpha;\beta,\nu)$ for every $\nu$.

The following two symmetries save us much work in the proofs for the classification of multiplicity-free skew characters. The first symmetry is $c(\lambda';\mu',\nu')=c(\lambda;\mu,\nu)$ which we will call conjugation symmetry. Further we have the rotation symmetry $[(\lambda/\mu)^\circ]=[\lambda/\mu]$.

We say that a skew diagram $\mathcal{D}$ decomposes into the disconnected skew diagrams $\mathcal{A}$ and $\mathcal{B}$ if no box of $\mathcal{A}$ (viewed as boxes in $\mathcal{D}$) is in the same row or column as a box of $\mathcal{B}$.

We will call a skew diagram a \textbf{proper skew diagram} if it is neither a partition nor a partition rotated by $180^\circ$. It is known that $[\lambda/\mu]=[\nu]$ for some irreducible character $[\nu]$ when $\lambda/\mu$ is a partition or a partition rotated by $180^\circ$, but if $\lambda/\mu$ is a proper skew diagram then there occur at least two different irreducible characters in the decomposition of $[\lambda/\mu]$ into irreducible characters (see \cite{Bes}).

Stembridge classified in \cite{Ste} the multiplicity-free products of Schur functions $s_\mu\cdot s_\nu$. We call a partition $\lambda$ a $k$-line rectangle if $\lambda=(i^k)$ or $\lambda=(k^i)$ for some $i\geq1$. A partition is called a near-rectangle if it is possible to obtain a rectangle by deleting a single column or row. A partition $\lambda$ is a fat hook if $\lambda=(\lambda_1^{l_1},\lambda_2^{l_2})$ for some $l_1,l_2\geq1, \lambda_1>\lambda_2\geq1$. So near-rectangles are special fat hooks. If we exclude the trivial case that one partition is empty the classification of multiplicity-free products of Schur functions is as follows (Theorem 2.1 in \cite{Ste}):

\begin{Sa}
\label{Sa:Mainstem}
Let $\mu, \nu$ be partitions. The product $s_\mu\cdot s_\nu$ of two Schur functions is multiplicity-free if and only if up to exchange of $\mu$ and $\nu$ one of the following conditions holds:
\begin{enumerate}
		\item $\mu$ is a one-line rectangle
		\item $\mu$ is a two-line rectangle and $\nu$ is a fat hook
		\item $\mu$ is a rectangle and $\nu$ a near-rectangle
		\item $\mu$ and $\nu$ are rectangles
\end{enumerate}
\end{Sa}

This theorem will be important in our proofs for the classification of  mul\-ti\-plicity-free skew characters.

\section{Multiplicity-free skew characters}
\label{sec:skew}
In this section we will introduce a partial order on the set of skew characters which we use to classify the multiplicity-free skew characters.

For this we need a generalization of the theorem which Stembridge used in \cite{Ste} to classify multiplicity-free products of Schur functions. Stembridge proved the following theorem for the special cases where $a=b$ or $b=0$.

\begin{Sa}
\label{Sa:schieb}
	Let $\lambda, \mu, \nu $ be partitions and $a\geq b\geq 0$ integers. Then:
	\[c(\lambda;\mu,\nu) \leq c(\lambda+(1^a);\mu+(1^b),\nu+(1^{a-b}))\]
	as well as
	\[c(\lambda;\mu,\nu) \leq c(\lambda \cup (a);\mu \cup (b),\nu \cup (a-b))\]
\end{Sa}
{\itshape Proof: } We will show how to obtain from a LR-tableau of shape $\lambda/\mu$ with content $\nu$ a LR-tableau of shape $\left(\lambda+(1^a)\right)/\left(\mu+(1^b)\right)$ with content $\nu+(1^{a-b})$. Then we will argue that the new tableaux are all different and from this follows the first inequality. The second follows from the first by conjugation symmetry. To shorten the proof we will assume that the boxes of $\lambda$ which usually should be deleted to get $\lambda/\mu$ are instead filled with zeros.

Let $\alpha$ be a LR-tableau of shape $\lambda/\mu$ with content $\nu$ and let $c=a-b$. Into each of the first $b$ rows we insert a zero such that the rows are semistandard, i.e. weakly increasing from left to right. For $1\leq i \leq c$ we insert a $i$ into the row $b+i$ such that the rows are again semistandard. This tableau we will call $\bar\alpha$ and we have to show, that $\bar\alpha$ is a LR-tableau.

The reverse row word of $\bar\alpha$ is still a lattice word, because before every new entry $i>1$ in the reverse row word there is also a new entry $i-1$.

Clearly $\bar\alpha$ satisfies the semistandard conditions for the rows. The boxes labelled $0$ in $\bar\alpha$ form the partition $\bar\mu=\mu+(1^b)$.

{\bfseries Example:} For $\lambda=(7^2,5^3,2,1),\mu=(6,3^2),\nu=(6,5^2,4)$ and $a=6,b=2$ we have the following start tableaux (1) and end tableaux (2).  (Here $c(\lambda;\mu,\nu)=2\leq3=c(\lambda +(1^6);\mu+(1^2),\nu+(1^4))$.)
\newcommand{\redzero}{\mathbf{\textcolor{red}{0}}}
\newcommand{\redone}{\mathbf{\textcolor{red}{1}}}
\newcommand{\redtwo}{\mathbf{\textcolor{red}{2}}}
\newcommand{\redthree}{\mathbf{\textcolor{red}{3}}}
\newcommand{\redfour}{\mathbf{\textcolor{red}{4}}}

{
\centering \footnotesize
\[1.: \young(::::::1,:::1112,:::22,11333,22444,33,4)=\young(0000001,0001112,00022,11333,22444,33,4) \quad \left| \quad \young(0000001,0001112,00022,11233,23344,34,4)=\young(::::::1,:::1112,:::22,11233,23344,34,4)\right.\]

\[2.: \young(:::::::1,::::1112,:::122,112333,223444,334,4)=\young(000000\redzero1,000\redzero1112,000\redone22,11\redtwo333,22\redthree444,33\redfour,4) \quad \left| \quad \young(000000\redzero1,000\redzero1112,000\redone22,112\redtwo33,233\redthree44,34\redfour,4)=\young(:::::::1,::::1112,:::122,112233,233344,344,4)\right.\]

}

\vspace{0.5cm}

To show that $\bar\alpha$ is a LR-tableau we now have to show that the semistandard conditions for the columns of $\bar\alpha$ are also satisfied:

For the following let $\kappa_j(i)$ be the number of entries smaller or equal to $i$ in the $j$-th row (we also count the zeros). Clearly a tableau with rows obeying the semistandard conditions is semistandard if we have for each $i\geq1$ and $j\geq2$: $\kappa_{j-1}(i-1)\geq\kappa_j(i)$. If $\alpha$ is semistandard and therefore meets the conditions for the $\kappa$ then also $\bar\alpha$ will meet the conditions, because for some given $i$ and $j$ we either change neither $\kappa_{j-1}(i-1)$ nor $\kappa_j(i)$ or only $\kappa_{j-1}(i-1)$ changes by $+1$ (in the case $j-1=a$ and $i-1\geq a-b$) or both will be changed by $+1$ by going from $\alpha$ to $\bar\alpha$. Thus $\bar\alpha$ is semistandard if $\alpha$ is.

Hence we only have to show that different LR-tableaux $\alpha,\beta$ of shape $\lambda/\mu$ with content $\nu$ give different $\bar\alpha,\bar\beta$. For this we notice that the semistandard condition together with the information of the content of each row determines the filling of a LR-tableau, and when two LR-tableaux are different there are two rows with different content. But if there are rows with different content in $\alpha$ and $\beta$ and we add to each row in $\alpha$ and $\beta$ the same entry, then there are also rows with different content in $\bar\alpha$ and $\bar\beta$, and so $\bar\alpha$ and $\bar\beta$ are different LR-tableaux. \qed

\begin{Bem} \normalfont
This theorem will play the important role in our proof of the classification of multiplicity-free skew characters. If we start with a skew diagram $\mathcal{D}$ whose character is not mul\-ti\-plicity-free, then all diagrams we can get from $\mathcal{D}$ by adding skew rows ($(a)/(b)$) and/or skew columns ($(1^a)/(1^b)$) have non multiplicity-free characters. On the other hand if we obtain a diagram with multiplicity-free character by adding skew rows and/or columns to some diagram $\mathcal{D}$ the theorem tells us that $[\mathcal{D}]$ is multiplicity-free. In our proofs we will insert skew rows/columns until we get a skew diagram containing two disconnected diagrams so that we can use Stembridge's classification of multiplicity-free products of Schur functions (see Theorem~\ref{Sa:Mainstem}).
\end{Bem}

The following is a well known result, which we will also need for the proofs of mul\-ti\-plicity-free skew characters.

\begin{Le}
\label{Sa:druck}
	Let $\lambda=(\lambda_1^{k_1},\ldots,\lambda_j^{k_j}), \mu=(\mu_1,\ldots,\mu_l), \nu$ be partitions.
	\begin{enumerate}
		\item If $l \le k_i$ for some $0\leq i \leq j$ then for all $n \geq 0$:
		\[c(\lambda;\mu,\nu) = c(\lambda \cup (\lambda_i^n);\mu ,\nu \cup (\lambda_i^n))\]
		\item If  $\mu_1 \le \lambda_i-\lambda_{i+1}$ (as usual $\lambda_{j+1}=0$) for some $0\leq i \leq j$ then let $r_i=\sum_{a=1}^i k_a$ and for all $n \geq 0$:
		\[c(\lambda;\mu,\nu) = c(\lambda + (n^{r_i});\mu,\nu + (n^{r_i}))\]
	\end{enumerate}
\end{Le}
{\itshape Proof:} The statements are equivalent by conjugation symmetry, so we will prove only the first. Looking at a  LR-filling of $\bar \lambda= \lambda \cup (\lambda_i^n)$ with content $\mu$, we see that the additional rows must be empty, because we can fill at most $l$ boxes in one column. So the first statement follows. \qed
\begin{Bsp}\normalfont
If we have $\lambda=(4,3^3,1),\mu=(3^3),\nu=(4,1),i=2,n=1$ the corresponding LR-tableaux are:
\begin{align}
  \textnormal{shape } \lambda/\nu \textnormal{, content } \mu:\young(\hfil\hfil\hfil\hfil,\hfil11,122,233,3)&&\textnormal{shape } \bar\lambda/\bar\nu \textnormal{, content } \mu: \young(\hfil\hfil\hfil\hfil,\hfil\hfil\hfil,\hfil11,122,233,3)\nonumber
\end{align}
We see that the second row in $\bar\lambda$ must remain empty in any LR-tableau with filling $\mu$ because of the LR-conditions.
\end{Bsp}

We are now ready to give the classification of multiplicity-free skew characters and prove the classification.

For this we assume in the following that for partitions $\lambda=(\lambda_1,\ldots,\lambda_k),\alpha=(\alpha_1,\ldots,\alpha_l)$ the following conditions holds true $\alpha \subseteq\lambda,\alpha_1<\lambda_1, l<k$. These are no restrictions. The conditions merely allow us not to worry about nonexistent rows or columns and so keep the classification simple (a shorter classification is given in Theorem~\ref{Sa:Mainskewshort}). Furthermore if we write $\lambda=(\lambda_1^{k_1},\lambda_2^{k_2},\ldots\lambda_i^{k_i})$ then we assume, that we have $\lambda_h>\lambda_{h+1}$ for all $1\leq h \leq i-1$. Thus in the theorem $i$ will be the number of distinct part sizes of $\lambda$ and $j$ the number for $\alpha$.

\begin{Sa}
\label{Sa:Mainskew} Let $\lambda=(\lambda_1^{k_1},\lambda_2^{k_2},\ldots\lambda_i^{k_i}),\alpha=(\alpha_1^{l_1},\dots,\alpha_j^{l_j})$ be partitions and $k=\sum_{a=1}^ik_a, l=\sum_{a=1}^jl_a$.

Then: $[\lambda/\alpha]$ is multiplicity-free if and only if one of the following conditions holds:
	\begin{enumerate}
	\item $\lambda/\alpha$ decomposes into two disconnected skew diagrams $\mu$ and $\nu$ for which up to rotation by $180^\circ$ and/or exchanging $\mu$ with $\nu$ one of the following conditions holds:
		\begin{enumerate}
		\item $\mu$ is a one-line rectangle and $\nu$ is a partition
		\item $\mu$ is a two-line rectangle and $\nu$ is a fat hook
		\item $\mu$ is a rectangle and $\nu$ a near-rectangle
		\item $\mu$ and $\nu$ are rectangles
		\end{enumerate}
	\item $\lambda/\alpha$ is a connected skew diagrams and one of the following conditions holds:
		\begin{enumerate}
		\item $i=1$
		\item $j=1$ and one of the following conditions holds:
			\begin{enumerate}
			\item $\alpha_1=1$ or $l_1=1$
			\item $\lambda_1=1+\alpha_1$ or $k=1+l$
			\item $i=2$
			\item $i=3$ and one of the following conditions holds:
				\begin{enumerate}
				\item $\alpha_1=2$ or $l_1=2$
				\item $k_1=1$ or $\lambda_3=1$
				\item $k_2=1$ or $\lambda_2=1+\lambda_3$
				\item $k_3=1$ or $\lambda_1=1+\lambda_2$
				\item $k=2+l$ or $\lambda_1=2+\alpha_1$
				\end{enumerate}
			\end{enumerate}
		\item $i=2$ and one of the following conditions holds:
			\begin{enumerate}
			\item $\lambda_1=1+\lambda_2$ or $k_2=1$
			\item $\lambda_2=1$	or $k_1=1$
			\end{enumerate}

		\item $i=2$ and $j=2$ and one of the following conditions holds:
			\begin{enumerate}
			\item $\lambda_1=1+\alpha_1$ or $k=1+l$
			\item $\lambda_1=2+\lambda_2$ or $k_2=2$
			\item $\lambda_2=2$ or $k_1=2$
			\item $\alpha_1=1+\alpha_2$ or $l_2=1$
			\item $\alpha_2=1$ or $l_1=1$
			\end{enumerate}

		\end{enumerate}
	\end{enumerate}
\end{Sa}

\textbf{In the following we use the notation of Theorem~\ref{Sa:Mainskew}}.

\begin{pro}
\label{Sa:3teile}
$[\lambda/\alpha]$ is not multiplicity-free if $\lambda/\alpha$ decomposes into three (or more) disconnected skew diagrams.
\end{pro}
{\itshape Proof:} This is equivalent to the product of 3 Schur functions. Stembridge showed in \cite{Ste} that this is not multiplicity-free, by using $c((3,2,1);(2,1),(2,1))=2$ and a special case of Theorem~\ref{Sa:schieb}. So in this case $[\lambda/\alpha]$ is not multiplicity-free. \qed

\begin{pro}
\label{Sa:2teileskew}
$[\lambda/\alpha]$ is not multiplicity-free if $\lambda/\alpha$ decomposes into two disconnected skew diagrams $\mu/\beta$ and $\nu/\gamma$ and at least one of them is a proper skew partition.
\end{pro}
{\itshape Proof:} We will show how to obtain $\lambda/\alpha$ by adding skew rows and/or columns to  $(3,2,1)/(2,1)$. Here $c((3,2,1);(2,1),(2,1))=2$, because we have the following LR-tableaux: {\small$\young(::1,:2,1),\young(::1,:1,2)$}. Theorem~\ref{Sa:schieb} then tells us that $[\lambda/\alpha]$ is not multiplicity-free.

We may assume that $\mu/\beta$ is a proper skew partition and that $\mu/\beta$ is to the right and above $\nu/\gamma$. We know that $\mu$ has at least two parts $\mu_1\not=\mu_2$, because otherwise $\mu/\beta$ would be a rotated partition. Furthermore let $\beta_1$ denote the largest part of $\beta$ and $\beta_2\not=\beta_1$  one part of $\beta$ (including the case $\beta_2=0$) such that the skew row $(\mu_2)/(\beta_2)$ is a part of $\mu/\beta$, for example we could choose the last row of $\mu/\beta$. With $\nu_1$ and $\gamma_1$ we refer to the biggest parts of $\nu$ and $\gamma$ respectively and we assume that $\nu_1\not=\gamma_1$.

We begin with the diagram $(3,2,1)/(2,1)$ and add $\gamma_1$ skew columns $(1^3)/(1^3)$ pushing the diagram to the right to
\[(3+\gamma_1,2+\gamma_1,1+\gamma_1)/(2+\gamma_1,1+\gamma_1,\gamma_1).\]
Adding $\nu_1-1-\gamma_1\geq0$ skew columns $(1^3)/(1^2)$ gives:
\[(3+\nu_1-1,2+\nu_1-1,1+\nu_1-1)/(2+\nu_1-1,1+\nu_1-1,\gamma_1)\]
Adding $\beta_2$ skew columns $(1^2)/(1^2)$ gives:
\[(2+\beta_2+\nu_1,1+\beta_2+\nu_1,\nu_1)/(1+\beta_2+\nu_1,\beta_2+\nu_1,\gamma_1)\]
Now we must distinguish the two cases $\mu_2<\beta_1$ and $\mu_2\geq\beta_1$:

For $\mu_2<\beta_1$ we add $\mu_2-1-\beta_2\geq0$ skew columns $(1^2)/(1^1)$ and get:
\[(1+\mu_2+\nu_1,\mu_2+\nu_1,\nu_1)/(\mu_2+\nu_1,\beta_2+\nu_1,\gamma_1)\]
Now we add $\beta_1-\mu_2>0$ skew columns $(1^1)/(1^1)$ and then $\mu_1-1-\beta_1\geq0$ skew columns $(1^1)/(1^0)$ to get:
\[(1+\beta_1+\mu_1-1-\beta_1+\nu_1,\mu_2+\nu_1,\nu_1)/(\beta_1+\nu_1,\beta_2+\nu_1,\gamma_1)=\]\[(\mu_1+\nu_1,\mu_2+\nu_1,\nu_1)/(\beta_1+\nu_1,\beta_2+\nu_1,\gamma_1)\]

For $\mu_2\geq\beta_1$ we add $\beta_1-\beta_2-1\geq0$ skew columns $(1^2)/(1^1)$ and get:
\[(1+\beta_1+\nu_1,\beta_1+\nu_1,\nu_1)/(\beta_1+\nu_1,\beta_2+\nu_1,\gamma_1)\]
Now we add $\mu_2-\beta_1$ skew columns $(1^2)/(1^0)$ and then $\mu_1-\mu_2-1$ skew columns $(1^1)/(1^0)$ to get again:
\[(\mu_1+\nu_1,\mu_2+\nu_1,\nu_1)/(\beta_1+\nu_1,\beta_2+\nu_1,\gamma_1)\]

In both cases we have the same diagram with three skew rows, and these skew rows also appear in $\lambda/\alpha$. If $\lambda/\alpha$ has more skew rows, we may add them also to the diagram to get $\lambda/\alpha$. Because we have only used the operations covered in Theorem~\ref{Sa:schieb} the theorem tells us that there is a $\mu$ with $c(\lambda;\alpha,\mu)\geq2$ and so $[\lambda/\alpha]$ is not multiplicity-free. \qed

\begin{Bem} \label{Rem}
\begin{enumerate}
 \item The way from $(3,2,1)/(2,1)$ to $\lambda/\alpha$ will not be shown again in the following proofs for characters with multiplicity.
\item Proposition~\ref{Sa:3teile} follows from Proposition~\ref{Sa:2teileskew} if we view two of the three skew diagrams as one proper skew diagram.
\item \label{Rem3}The proof of Proposition~\ref{Sa:2teileskew} also gives us a character $[\mu]$ in  $[\lambda/\alpha]$ with multiplicity at least 2.
\end{enumerate}
\end{Bem}

\begin{Sa}
\label{Sa:2teile}
Let $\lambda/\alpha$ decompose into two disconnected skew diagrams $\mu$ and $\nu$. Then $[\lambda/\alpha]$ is multiplicity-free if and only if for $\mu$ and $\nu$ one of the following conditions holds up to rotation by $180^\circ$ and/or exchanging $\mu$ with $\nu$:
		\begin{enumerate}
		\item $\mu$ is a one-line rectangle and $\nu$ is a partition
		\item $\mu$ is a two-line rectangle and $\nu$ is a fat hook
		\item $\mu$ is a rectangle and $\nu$ a near-rectangle
		\item $\mu$ and $\nu$ are rectangles
		\end{enumerate}
\end{Sa}
{\itshape Proof:} If $\mu$ or $\nu$ is a proper skew diagram then $[\lambda/\alpha]$ is not multiplicity-free by Proposition~\ref{Sa:2teileskew}. If both $\mu$ and $\nu$ are partitions or partitions rotated by $180^\circ$  $[\lambda/\alpha]$ is equivalent to the product $s_\mu \cdot s_\nu$ of Schur functions. So we may use the classification of Stembridge Theorem~\ref{Sa:Mainstem} and the theorem follows. \qed
\vspace{0.5cm}

We have now covered all cases in which $\lambda/\alpha$ decomposes into disconnected skew diagrams. For the following lemmas we will assume that $\lambda/\alpha$ is connected. We will continue to use the notation of Theorem~\ref{Sa:Mainskew}.

\begin{Le}
\label{Le:i1}
If $i=1$ then $[\lambda/\alpha]$ is multiplicity-free.
\end{Le}
{\itshape Proof:} In this case $\lambda/\alpha$ is a partition rotated by $180^\circ$, so $[\lambda/\alpha]$ is multiplicity-free (see \cite{Bes}). \qed

\begin{Le}
\label{Le:a1line}
Let $j=1$.

If $\alpha_1=1$ or $l_1=1$ then $[\lambda/\alpha]$ is multiplicity-free.

If $\lambda_1=1+\alpha_1$ or $k=1+l$ then $[\lambda/\alpha]$ is multiplicity-free.
\end{Le}
{\itshape Proof:} Let $\alpha_1=1$. In this case $\alpha$ is a one column rectangle. Inserting $k-l$ skew rows $(1)/(1)$, this means we push the first column downwards, we get a skew diagram $\mu/\beta$ containing two disconnected skew diagrams, one of which is a one column rectangle and the other is some partition (e.g.:{\tiny $\young(:\hfil\hfil\hfil\hfil,:\hfil\hfil\hfil\hfil,:\hfil\hfil,\hfil\hfil,\hfil)$} becomes {\tiny $\young(:\hfil\hfil\hfil\hfil,:\hfil\hfil\hfil\hfil,:\hfil\hfil,:\hfil,\hfil,\hfil)$}). By Theorem~\ref{Sa:2teile} $[\mu/\beta]$ is multiplicity-free and by Theorem~\ref{Sa:schieb} $[\lambda/\alpha]$ is also multiplicity-free. The case $l_1=1$ is equivalent to the case $\alpha_1=1$ by conjugation symmetry.

Let $\lambda_1=1+\alpha_1$. By inserting a skew column $(1^l)/(1^l)$ we push the rectangle (initially still connected to the other boxes) to the right and get a skew diagram which decomposes into two disconnected skew diagrams, one of which is a one column rectangle and the other is some partition (e.g.:{\tiny $\young(::::\hfil,::::\hfil,::::\hfil,\hfil\hfil\hfil\hfil\hfil,\hfil\hfil\hfil,\hfil\hfil)$} becomes {\tiny $\young(:::::\hfil,:::::\hfil,:::::\hfil,\hfil\hfil\hfil\hfil\hfil,\hfil\hfil\hfil,\hfil\hfil)$}). So $[\lambda/\alpha]$ is multiplicity-free. The case $k=1+l$ is equivalent to the case $\lambda_1=1+\alpha_1$ by conjugation symmetry.\qed

\begin{Le}
\label{Le:i2-1line}
Let $i=2$.

If $\lambda_1=1+\lambda_2$ or $k_2=1$ then $[\lambda/\alpha]$ is multiplicity-free.

If $\lambda_2=1$	or $k_1=1$ then $[\lambda/\alpha]$ is multiplicity-free.
\end{Le}
{\itshape Proof:} This is equivalent to the cases in Lemma~\ref{Le:a1line} by rotation symmetry. For example, rotating a skew diagram with $i=2$ and $\lambda_1=1+\lambda_2$  by $180^\circ$ yields a skew diagram with $j=1$ and $\alpha_1=1$ (e.g.:{\tiny $\young(::::\hfil\hfil,:::\hfil\hfil\hfil,::\hfil\hfil\hfil,\hfil\hfil\hfil\hfil\hfil)$} becomes {\tiny $\young(:\hfil\hfil\hfil\hfil\hfil,:\hfil\hfil\hfil,\hfil\hfil\hfil,\hfil\hfil)$}). \qed

\begin{Le}
\label{Le:j1i2}
Let $j=1$ and $i=2$.

Then $[\lambda/\alpha]$ is multiplicity-free.\end{Le}
{\itshape Proof:} We may assume that $\alpha_1=\lambda_2$: If we have $\alpha_1 < \lambda_2$, then let $\bar\lambda$ denote the partition one gets by removing $\lambda_2-\alpha_1$ columns $(1^k)$ of $\lambda$. Lemma~\ref{Sa:druck} tells us that $[\lambda/\alpha]$ is multiplicity-free if and only if $[\bar\lambda/\alpha]$ is multiplicity-free. If we have $\alpha_1>\lambda_2$ then we may add $\alpha_1-\lambda_2$ columns $(1^k)$ to $\lambda$ to get $\bar\lambda$ and Theorem~\ref{Sa:schieb} tells us that $[\lambda/\alpha]$ is multiplicity-free if $[\bar\lambda/\alpha]$ is. Hence let $\alpha_1=\lambda_2$. If we add rows $(\alpha_1)/(\alpha_1)$ to $\lambda/\alpha$ until we have two disconnected skew diagrams, then both disconnected skew diagrams are rectangles and the skew character is multiplicity-free by Theorem~\ref{Sa:2teile}. By Theorem~\ref{Sa:schieb} $[\lambda/\alpha]$ is multiplicity-free.\qed

\begin{Bem}  Reasoning similar to that justifying the claim that $\alpha_1=\lambda_2$ is not a loss of generality will not be repeated in the following proofs.
\end{Bem}

\begin{Le}
\label{Le:j1i3}
Let $j=1$ and $i=3$.

$[\lambda/\alpha]$ is multiplicity-free if and only if one of the following conditions holds:
\begin{enumerate}
				\item \label{l131}$\alpha_1=1$ or $l_1=1$
				\item \label{l132}$k=1+l$ or $\lambda_1=1+\alpha_1$
				\item \label{l133}$\alpha_1=2$ or $l_1=2$
				\item \label{l134}$k_1=1$ or $\lambda_3=1$
				\item \label{l135}$k_2=1$ or $\lambda_2=1+\lambda_3$
				\item \label{l136}$k_3=1$ or $\lambda_1=1+\lambda_2$			\item \label{l137}$k=2+l$ or $\lambda_1=2+\alpha_1$
\end{enumerate}
\end{Le}
{\itshape Proof:} If none of the mentioned conditions holds, then we get $\lambda/\alpha$ by adding skew rows and/or columns in a similiar way as in Proposition~\ref{Sa:2teileskew} to {$(6^2,4^2,2^2)/(3^3)=$}{ \tiny $\young(:::\hfil\hfil\hfil,:::\hfil\hfil\hfil,:::\hfil,\hfil\hfil\hfil\hfil,\hfil\hfil,\hfil\hfil)$}. Here $c((6^2,4^2,2^2);(3^3),(5,4,3,2,1))=2$ because we have the following LR-tableaux: {\scriptsize$\young(:::111,:::222,:::3,1134,23,45) \qquad \young(:::111,:::222,:::3,1134,24,35)$}. By Theorem~\ref{Sa:schieb}, $[\lambda/\alpha]$ is not multiplicity-free.

Now to the multiplicity-free cases:

\begin{itemize}
				\item \ref{l131} \& \ref{l132}:  The cases $\alpha_1=1$, $l_1=1$, $k=1+l$ and $\lambda_1=1+\alpha_1$ are covered by Lemma~\ref{Le:a1line}.
				\item \ref{l133}(b): $l_1=2$: We may assume that $k_1=2$. Adding skew columns $(1^2)/(1^2)$ pushes the first two rows to the right. If we push until we have two disconnected skew diagrams, then one is a two line rectangle and the other is a fat hook. The character is multiplicity-free by Theorem~\ref{Sa:2teile} and therefore $[\lambda/\alpha]$ is multiplicity-free by Theorem~\ref{Sa:schieb}.
				\item \ref{l133}(a): $\alpha_1=2$: This is equivalent to the case $l_1=2$ by conjugation symmetry.
				\item \ref{l134}(a): $k_1=1$: We may assume that $\alpha_1=\lambda_3$. Adding skew rows $(\alpha_1)/(\alpha_1)$ pushes a rectangle of width $\alpha_1$ to the bottom. If we push until we get two disconnected skew diagrams, then one is a rectangle and the other one is a near-rectangle. The character is multiplicity-free by Theorem~\ref{Sa:2teile} and therefore $[\lambda/\alpha]$ is multiplicity-free by Theorem~\ref{Sa:schieb}.
				\item \ref{l134}(b): $\lambda_3=1$: This is equivalent to the case $k_1=1$ by conjugation symmetry.
				\item \ref{l135}(a): $k_2=1$: We may assume that $\alpha_1=\lambda_3$. Pushing the lower rectangle of width $\alpha_1$ to the bottom yields two disconnected skew diagrams, one of which is a rectangle and the other is a near rectangle. The character is multiplicity-free by Theorem~\ref{Sa:2teile} and therefore $[\lambda/\alpha]$ is multiplicity-free by Theorem~\ref{Sa:schieb}.
				\item \ref{l135}(b): $\lambda_2=1+\lambda_3$: This is equivalent to the case  $k_2=1$ by conjugation symmetry.
				\item \ref{l136}(a): $k_3=1$: We may assume that $l=k_1$. Pushing the upper rectangle of height $k_1$ to the right yields a rectangle and a near-rectangle, so $[\lambda/\alpha]$ is multiplicity-free.
				\item \ref{l136}(b): $\lambda_1=1+\lambda_2$: This is equivalent to the case  $k_3=1$ by conjugation symmetry.
				\item \ref{l137}(a): $\lambda_1=2+\alpha_1$: We may assume that $l=k_1$. Pushing the upper rectangle of width $2$ to the right yields a fat hook and a two-line rectangle, so  $[\lambda/\alpha]$ is multiplicity-free.
				\item \ref{l137}(b): $k=2+l$: This is equivalent to the case  $\lambda_1=2+\alpha_1$ by conjugation symmetry.
				\qed
\end{itemize}

\begin{Le}
\label{Le:j2i2}
Let $j=2$ and $i=2$.

$[\lambda/\alpha]$ is multiplicity-free if and only if one of the following conditions holds:
\begin{enumerate}
			\item $\lambda_1=1+\lambda_2$ or $k_2=1$
			\item $\lambda_2=1$ or $k_1=1$
			\item $\lambda_1=2+\lambda_2$ or $k_2=2$
			\item $\lambda_1=1+\alpha_1$ or $k=1+l$
			\item $\alpha_1=1+\alpha_2$ or $l_2=1$
			\item $\alpha_2=1$ or $l_1=1$
			\item $\lambda_2=2$ or $k_1=2$
			\end{enumerate}
\end{Le}
{\itshape Proof:} This is equivalent to Lemma~\ref{Le:j1i3} by rotation symmetry. \qed

\begin{Le}
\label{Le:j1i4}
Let $j=1$ and $i\geq4$.

$[\lambda/\alpha]$ is multiplicity-free if and only if one of the following conditions holds:
\begin{enumerate}
		\item $\alpha_1=1$ or	$l_1=1$
		\item $\lambda_1=1+\alpha_1$ or $k=1+l$
\end{enumerate}
\end{Le}
{\itshape Proof:} If none of the mentioned conditions holds, then we get $\lambda/\alpha$ by adding skew rows and/or columns to {$(4,3,2,1)/(2^2)=$\tiny $\young(::\hfil\hfil,::\hfil,\hfil\hfil,\hfil)$}. Here $c((4,3,2,1);(2^2),(3,2,1))=2$ because we have the following LR-tableaux: {\scriptsize$\young(::11,::2,12,3) \qquad \young(::11,::2,13,2)$}. By Theorem~\ref{Sa:schieb}, $[\lambda/\alpha]$ is not multiplicity-free.

The multiplicity-free cases are those covered by Lemma~\ref{Le:a1line}.\qed

\begin{Le}
\label{Le:j3i2}
Let $j\geq3$ and $i=2$.

$[\lambda/\alpha]$ is multiplicity-free if and only if one of the following conditions holds:			\begin{enumerate}
			\item $\lambda_1=1+\lambda_2$ or $k_2=1$
			\item $\lambda_2=1$	or $k_1=1$
			\end{enumerate}
\end{Le}
{\itshape Proof:} This is equivalent to Lemma~\ref{Le:j1i4} by rotation symmetry. \qed

\begin{Le}
\label{Le:j2i3}
Let $j\geq2$ and $i\geq3$.

$[\lambda/\alpha]$ is not multiplicity-free.
\end{Le}
{\itshape Proof:} We obtain $\lambda/\alpha$ by adding skew rows and columns to {$(3,2,1)/(2,1)=$\tiny $\young(::\hfil,:\hfil,\hfil)$}. Here $c((3,2,1);(2,1),(2,1))=2$ as mentioned above. By Theorem~\ref{Sa:schieb}, $[\lambda/\alpha]$ is not multiplicity-free.\qed

\section{Schubert Calculus}
\label{sec:Schubert}
The cohomology ring $H^*(Gr(l,\mathbb{C}^n),\mathbb{Z})$ of the Grassmannian $Gr(l,\mathbb{C}^n)$ of $l$-di\-men\-sio\-nal subspaces of $\mathbb{C}^n$ has an additive basis of Schubert classes $\sigma_\lambda$. These are indexed by partitions $\lambda \subseteq (k^l)$ where $k=n-l$. Looking at the product of Schubert classes, we get the following formula:
\[\sigma_\mu\cdot\sigma_\nu=\sum_{\lambda\subseteq(k^l)}c(\lambda;\mu,\nu)\sigma_\lambda\]

The following lemma will provide us with another symmetry of the LR-coeffi\-cients with which we will show how to obtain the decomposition of a product of Schubert classes from the decomposition of its associated skew character (and vice versa).

\begin{Le}
\label{Sa:skewschub}
Let $\lambda,\mu,\nu$ be partitions with $\nu\subseteq\lambda$ and $\lambda\subset (k^l) \not= \lambda$ for some positive integers $k,l$. Then:
\[c(\lambda;\mu,\nu)=c\left((k^l/\nu)^\circ;\mu,(k^l/\lambda)^\circ\right)\]
\end{Le}
{\itshape Proof:} Let $\lambda^{-1}=(k^l)/\lambda$, $\nu^{-1}=(k^l)/\nu$. If we put $\lambda/\nu$ into a rectangle with $l$ rows and $k$ columns then $\lambda^{-1}$ is a partition rotated by $180^\circ$ in the lower right corner with the shape of $\lambda^{-1}$ running along the shape of $\lambda$ and correspondingly for $\nu^{-1}$.

If we look at the difference $\nu^{-1}/\lambda^{-1}$ we get back the skew diagram $\lambda/\nu$. If we now rotate $\lambda^{-1}$ and $\nu^{-1}$ to get partitions and look at ${\nu^{-1}}^\circ/{\lambda^{-1}}^\circ$ we get the skew diagram $\lambda/\nu$ rotated by $180^\circ$.

{\bfseries Example:} For $\lambda/\nu=(4,2)/(1)=$ {\footnotesize $\young(:\hfil\hfil\hfil,\hfil\hfil)$} and $(k^l)=(4^3)$ we have:\linebreak
$\nu^{-1}=${\footnotesize $\young(:\hfil\hfil\hfil,\hfil\hfil\hfil\hfil,\hfil\hfil\hfil\hfil)$} and $\lambda^{-1}=${\footnotesize $\young(::\hfil\hfil,\hfil\hfil\hfil\hfil)$} and ${\nu^{-1}}^\circ/{\lambda^{-1}}^\circ=$ {\footnotesize $\young(::\hfil\hfil,\hfil\hfil\hfil)$} which is $\lambda/\nu$ rotated by $180^\circ$.

Therefore we have $[\lambda/\nu]=[{\nu^{-1}}^\circ/{\lambda^{-1}}^\circ]$ (see \cite{Bes}) and so for every partition $\mu$ the coefficient of $[\mu]$ in the decomposition of both skew characters is the same, giving:
\[c(\lambda;\mu,\nu)=c\left((k^l/\nu)^\circ;\mu,(k^l/\lambda)^\circ\right) \qed\]

If we assume in Theorem~\ref{Sa:Mainskewschub} that $\mu \subseteq {\left((k^l)/\nu\right)}^\circ \not= \mu$ we only exclude the trivial products of Schubert classes $\sigma_\mu\cdot\sigma_\nu=0$ and $\sigma_\mu\cdot\sigma_\nu=\sigma_{(k^l)}$.

\begin{Sa}
\label{Sa:Mainskewschub}
Let $\mu,\nu,\lambda$ be partitions with $\mu,\nu\subseteq (k^l)$ for some positive integers $k,l$. Let $\mu \subseteq {\left((k^l)/\nu\right)}^\circ \not= \mu$.

Then: The coefficient of $\sigma_\lambda$ in the decomposition of the product $\sigma_\mu\cdot\sigma_\nu$ in $H^*(Gr(l,\mathbb{C}^n),\mathbb{Z})$ with $n=k+l$ equals the coefficient of $[((k^l)/\lambda)^\circ]$ in the decomposition of the skew character $[((k^l)/\nu)^\circ/\mu]$.
\end{Sa}
{\itshape Proof:}
We start with $[((k^l)/\nu)^\circ/\mu]$. The decomposition into irreducible characters is as follows:

\[ [((k^l)/\nu)^\circ/\mu] =\sum_{\alpha} c(((k^l)/\nu)^\circ;\mu,\alpha) [\alpha]\]

From the LR-Rule we see that $c(((k^l)/\nu)^\circ;\mu,\alpha)=0$ if $\alpha\not\subseteq((k^l)/\nu)^\circ$. Therefore it is sufficient to sum over $\alpha\subseteq((k^l)/\nu)^\circ$ or more generously over $\alpha\subseteq (k^l)$:

\[ [(k^l/\nu)^\circ/\mu] =\sum_{\alpha \subseteq (k^l)} c((k^l/\nu)^\circ;\mu,\alpha) [\alpha]\]

Using Lemma~\ref{Sa:skewschub} we get:

\[ [(k^l/\nu)^\circ/\mu] =\sum_{\alpha \subseteq (k^l)} c((k^l/\alpha)^\circ;\mu,\nu) [\alpha]\]

If we now set $(k^l/\alpha)^\circ=\lambda$ we get:

\[ [(k^l/\nu)^\circ/\mu] =\sum_{\lambda \subseteq (k^l)} c(\lambda;\mu,\nu) [(k^l/\lambda)^\circ]\]

Comparing this with the decomposition of the product of Schubert classes

\[ \sigma_\mu\cdot\sigma_\nu=\sum_{\lambda\subseteq (k^l)} c(\lambda;\mu,\nu) \sigma_\lambda \]

finishes the proof.\qed \newline

Thus instead of calculating the product of Schubert classes directly, we can decompose the skew character $[(k^l/\nu)^\circ/\mu]$ and then replace all characters $[\lambda]$ by Schubert classes $\sigma_{\left((k^l)/\lambda\right)^\circ}$. With this one can also compute products of Schur functions $s_\mu\cdot s_\nu$ by setting $k=\mu_1+\nu_1,l=l(\mu)+l(\nu)$. Because we have $c(\lambda;\mu,\nu)=0$ for $\lambda \not \subseteq (k^l)$ the sum $\sum_{\lambda\subseteq (k^l)} \sigma_\lambda$ becomes $\sum_\lambda \sigma_\lambda$. But this is just the sum for the product of Schur functions.

To the product $\sigma_\mu\cdot\sigma_\nu$ we associate a skew diagram in the following way. Remove from the rectangle $(k^l)$ the partition $\mu$ as usual and from the lower right corner remove $\nu^\circ$ to get a skew diagram. The associated \textbf{basic skew diagram} to the product $\sigma_\mu\cdot\sigma_\nu$ is then obtained by removing all empty rows and columns. The associated basic skew diagram is well defined unless $\mu$ und $\nu^\circ$ intersect in $(k^l)$, but then we have the trivial product: $\sigma_\mu\cdot\sigma_\nu=0$.

The associated basic skew diagram to an arbitrary skew diagram $\lambda/\mu$ is obtained by removing all empty rows and columns from $\lambda/\mu$.

We found the following idea of looking at the inner and outer lattice paths in \cite{Tho}.

For a basic skew diagram $\lambda/\mu$ (which is a proper skew diagram) we define two lattice paths from the lower left corner to the upper right corner. The outer lattice path starts to the right, follows the shape of $\lambda$ and ends upwards in the corner, while the inner lattice path starts upwards, follows the shape of $\mu$ and ends with a segment to the right. With $s_{in}$ we refer to the length of the shortest straight segment of the inner lattice path while $s_{out}$ is the length of the shortest straight segment of the outer lattice path.

{\bfseries Examples:}
\begin{enumerate}
\item The skew diagram {\tiny $\young(::::\hfil\hfil\hfil,::::\hfil\hfil\hfil,::\hfil\hfil\hfil\hfil\hfil,::\hfil\hfil\hfil,\hfil\hfil\hfil\hfil\hfil)$} has the inner lattice path lengths: 1,2,2,2,2,3 and therefore $s_{in}=1$ and the outer lattice path lengths: $5,2,2,3$ and therefore $s_{out}=2$.

\item The skew diagram {\tiny $\young(:::\hfil\hfil\hfil,:::\hfil\hfil\hfil,:::\hfil,:::\hfil,\hfil\hfil,\hfil\hfil,\hfil\hfil)$}  is not a basic skew diagram. \newline Its associated basic skew diagram {\tiny $\young(::\hfil\hfil\hfil,::\hfil\hfil\hfil,::\hfil,::\hfil,\hfil\hfil,\hfil\hfil,\hfil\hfil)$} has the inner lattice path lengths: $3,2,4,3$ (hence $s_{in}=2$) and the outer lattice path lengths: $2,3,1,2,2,2 $ (hence $s_{out}=1$).
\end{enumerate}

In this notation the classification of multiplicity-free skew characters looks like Stembridge's classification of multiplicity-free products of Schur functions (or equivalently the classification of multiplicity-free skew characters, whose diagram decomposes into two partitions). The following classification follows directly from Theorem~\ref{Sa:Mainskew}.

\begin{Sa}
\label{Sa:Mainskewshort}
Let $\alpha/\beta$ be a proper skew partition. Then $[\alpha/\beta]$ is multiplicity-free if and only if one of the following conditions holds for its associated basic skew diagram $\lambda/\mu$ up to exchanging $\mu$ and $s_{in}$ for $\nu$ and $s_{out}$ where $\nu={\left((\lambda_1^{l(\lambda)})/\lambda\right)}^\circ$:
\begin{enumerate}
\item $\mu$ is a rectangle and $s_{in}=1$
\item $\mu$ is a rectangle and $s_{in}=2$, $\nu$ is a fat hook
\item $\mu$ is a rectangle, $\nu$ is a fat hook and $s_{out}=1$
\item $\mu$ and $\nu$ are rectangles.
\end{enumerate}
\end{Sa}

There are three trivial cases for multiplicity-free products of Schubert classes:
\begin{enumerate}
\item If $\mu$ and $\nu^\circ$ intersect if $\mu$ is placed as usual in the upper left corner and $\nu^\circ$ is placed in the lower right corner of $(k^l)$, then $\sigma_\mu\cdot\sigma_\nu=0$.
\item If the associated basic skew diagram is empty, then $\sigma_\mu\cdot\sigma_\nu=\sigma_{(k^l)}$.
\item If the associated basic skew diagram is not a proper skew diagram, then $\sigma_\mu\cdot\sigma_\nu=\sigma_\lambda$ for some partition $\lambda$. (See Theorem~\ref{Sa:Mainskewschub} on how to obtain $\lambda$.)
\end{enumerate}

Thomas and Yong classified in \cite{Tho} the mul\-ti\-pli\-ci\-ty-free products of Schubert classes. Using Theorem~\ref{Sa:Mainskewshort} and Theorem~\ref{Sa:Mainskewschub} we get Thomas and Yong's classification of mul\-ti\-pli\-ci\-ty-free products of Schubert classes:

\begin{Sa}
\label{Sa:Mainschubert}
Let $\sigma_\alpha\cdot\sigma_\beta$ be a non-trivial (see above) product of two Schubert classes in $H^*(Gr(l,\mathbb{C}^n),\mathbb{Z})$.
Then  $\sigma_\alpha\cdot\sigma_\beta$ is multiplicity-free if and only if one of the following conditions holds for its associated basic skew diagram $\lambda/\mu$ up to exchanging $\mu$ and $s_{in}$ for $\nu$ and $s_{out}$ where $\nu={\left((\lambda_1^{l(\lambda)})/\lambda\right)}^\circ$:
\begin{enumerate}
\item $\mu$ is a rectangle and $s_{in}=1$
\item $\mu$ is a rectangle and $s_{in}=2$, $\nu$ is a fat hook
\item $\mu$ is a rectangle, $\nu$ is a fat hook and $s_{out}=1$
\item $\mu$ and $\nu$ are rectangles.
\end{enumerate}
\end{Sa}

\vspace{1cm}

{\bfseries Acknowledgement} This paper is based on the research I did for my diploma thesis supervised by Prof. Christine Bessenrodt. I am very grateful to Christine Bessenrodt for introducing me to this interesting field of algebra, supervising my research and helping me to write this paper.

Furthermore, thanks go to John Stembridge for drawing our attention to \cite{Tho}.

\end{document}